\documentclass[12pt]{amsart}
\usepackage{amsfonts,amssymb,amscd,xy,mathtools,tikz-cd}
\usepackage{setspace}
\usepackage[leqno]{amsmath}
\usepackage{mathrsfs}
\usepackage{amssymb}
\usepackage{amsmath}
\usepackage{amsxtra}
\usepackage{amsthm}
\usepackage{hyperref}
\usepackage{cleveref}
\newcommand{\rkk}{R_{\e k^{\le\prime}\!/k}}

\def\br{{\rm{Br}}\e}

\usepackage{rotating}
\usepackage{anyfontsize}
\fontsize{1}{2}

\newcommand{\Hom}{{\mathrm{Hom}}}

\DeclareMathAlphabet{\mathbbmsl}{U}{bbm}{m}{sl}

\newcommand{\car}{{\rm char}\,}

\def\sm{\lbe{\rm sm}}

\def\G{\mathbb{G}}

\def\g{\varGamma}

\def\s{\mathscr }

\newcommand{\into}{\hookrightarrow}
\newcommand{\onto}{\twoheadrightarrow}
\newcommand{\isoto}{\overset{\!\sim}{\to}}

\newcommand{\kp}{k^{\le\prime}}

\usepackage{epsfig,amssymb}
\usepackage{bm} 
\usepackage{verbatim} 

\definecolor{labelkey}{rgb}{1,0,0}

\usepackage{verbatim} 

\makeatletter

\DeclareMathAlphabet{\mathcalligra}{T1}{calligra}{m}{n}

\numberwithin{equation}{section}

\DeclareFontEncoding{OT2}{}{} 
\newcommand{\textcyr}[1]{%
{\fontencoding{OT2}\fontfamily{wncyr}\fontseries{m}\fontshape{n}
		\selectfont #1}}
\newcommand{\Sha}{{\!\be\lbe\mbox{\textcyr{Sh}}}}

\newcommand{\Bha}{{\!\be\lbe\mbox{\textcyr{B}}}}

\newcommand{\sh}{\kern -.4em\phantom{a}^{\mathbf{\sim}}}

\newcommand{\lra}{\longrightarrow}

\newcommand{\fl}{\le{\rm fl}}

\def\be{\kern -.1em}
\def\le{\kern 0.03em}
\def\lle{\kern 0.015em}
\def\lbe{\kern -.05em}

\newcommand{\A}{\mathbb A}

\newcommand{\Z}{{\mathbb Z}}
\newcommand{\Q}{{\mathbb Q}}

\newcommand{\R}{{\mathbb R}}

\newcommand{\spec}{\mathrm{ Spec}\,}

\newcommand{\krn}{\mathrm{Ker}\,}
\newcommand{\img}{\mathrm{Im}\e\le}

\newcommand{\cok}{\mathrm{Coker}\,}

\def\e{\kern 0.08em}

\newcommand{\alb}{{\rm alb}}

\usepackage{tabularx}
\usepackage{epsfig,amssymb}
\usepackage{bm} 
\usepackage{verbatim} 
\usepackage{mathrsfs}

\def\pic{{\rm{Pic}}\e}

\usepackage{hyperref}

\definecolor{labelkey}{rgb}{1,0,0}

\makeatletter

\DeclareMathAlphabet{\mathcalligra}{T1}{calligra}{m}{n}

\numberwithin{equation}{section}

\newcommand{\extk}{{\mathscr E\be x\le t}_{k_{\lbe\lbe\fl}^{\lle\sim}}}

\def\be{\kern -.1em}
\def\le{\kern 0.055em}
\def\lle{\kern 0.015em}
\def\lbe{\kern -.03em}

\newcommand*{\defeq}{\mathrel{\vcenter{\baselineskip0.5ex \lineskiplimit0pt
			\hbox{\scriptsize.}\hbox{\scriptsize.}}}%
	=}

\newcommand{\ksep}{k^{\le\rm s}}

\newtheorem{lemma}{Lemma}[section]
\newtheorem{theorem}[lemma]{Theorem}
\newtheorem{proposition-definition}[lemma]{Proposition-Definition}
\newtheorem{theorem-definition}[lemma]{Theorem\le-Definition}

\theoremstyle{definition}

\newtheorem{lemma-definition}[lemma]{Lemma-Definition}
\newtheorem{question}[lemma]{Question}

\theoremstyle{remark}
\newtheorem{remark}[lemma]{Remark}
\newtheorem{remarks}[lemma]{Remarks}
\newtheorem{example}[lemma]{Example}

\begin{document}

\input xy     
\xyoption{all} 

\title[On the new theory of $1\e$-\e motives]{On the new theory of $1\e$-\e motives}

\author{Cristian D. Gonz\'alez\e-\le Avil\'es}
\address{Departamento de Matem\'aticas, Universidad de La Serena, Cisternas 1200, La Serena 1700000, Chile}
\email{cgonzalez@userena.cl}
\thanks{The author is partially supported by Fondecyt grant 1200118.}

\maketitle

\begin{abstract} This is a letter (not intended for publication in a regular journal) written in response to two referees of my paper \cite{ga23}. In it, I discuss possible applications of the new theory of $1$-motives introduced in \cite{ga23}. One of these is a new approach to the BSD conjecture for abelian varieties over global fields (developed here only partially) which builds on Bloch's well-known volume-theoretic interpretation of the conjecture in the number field case \cite{blo80}.
\end{abstract}

\section{Introduction}

\topmargin -1cm

Dear Referees X and Y:

\medskip

Many thanks for your helpful comments, which are reproduced below. After reading your comments, it became clear to me that I did not explain well the significance and/or possible applications of the ideas (and methods) developed in \cite{ga23}. I try to remedy that here. 

\bigskip

\textbf{Referee X}: ``In the paper \cite{ros18} (a significant input into the submitted paper \cite{ga23}), Rosengarten established a generalization of Tate duality over local and global fields (main point being the imperfect function field cases) with coefficients in an arbitrary commutative {\it affine} group scheme of finite type.  One of the main insights was to work with ``dual" in the sense of the Hom-functor into $\G_m$, which even on $\G_{a}$ is non-trivial since that has lots of sections over non-reduced artin local rings.

The submitted paper carries out a natural generalization over local fields:  dropping the affineness hypothesis (in the paper, “algebraic group” means “group scheme of finite type”, but the author doesn't clearly say that).  Since duality for abelian varieties doesn't rest on $\Hom(-,\G_m)$ but rather ${\rm Ext}^1(-,\G_m)$, to formulate a suitable notion of ``dual" beyond the affine case requires going beyond the $k$-group setting and instead working with a ``1-motive" notion (introduced long ago by Deligne for some special settings with tori and abelian varieties).

The submitted paper makes a super-general ``1-motive" definition (Section 3) over any field that allows including all commutative group schemes of finite type.  (Note that if $G$ is a connected commutative group scheme of finite type over a field $k$ of $\car p > 0$, it is an extension of an abelian variety by a possibly non-smooth affine group scheme of finite type.)  This includes a notion of dual (Def. 3.6), and the main results are perfect pairings in hypercohomology of such $1$-motives over non-archimedean local fields, as stated in Theorem 0.1.

Given the complexities in the definition of ``1-motive over $k$", it is natural to wonder (for the sake of suitability for ANT):  is this generalization for its own sake, or for some other purpose?

For example, there was arithmetic motivation behind \cite{ros18}, namely in separate work of his to prove a ``Tamagawa number formula" for adelic volumes of general smooth connected commutative affine or pseudo-reductive groups (duality theory in the general commutative case was a key tool for handling the non-commutative pseudo-reductive case); Rosengarten also gave counterexamples to show his restrictions were basically optimal.

The submitted paper gives no indication of being for a purpose beyond generalization for its own sake, and the proofs in Section 4 for the main result(s) seem to [be] merging the results of \cite{ros18} and homological arguments to get results in the wider ``1-motive duality" notion.  It is unclear whether Section 4 has sufficiently significant new ideas: is the main impact from this paper expanding generality by finding the right definition with which to merge prior results with homological formalism, or something more?  The submitted paper doesn't include an application of its local duality results, as may be desired for the paper to be at the level for ANT."

\bigskip

\textbf{Referee Y}: ``The author proves new duality theorems for the fppf cohomology of commutative group schemes over local fields, the main case of interest being that of positive characteristic. Here major technical difficulties arise from the fact that the group schemes under consideration may not be smooth, and even when they are they may not be written as extensions of an abelian variety by a smooth affine group scheme as in Chevalley's theorem. The author overcomes this problem by considering a certain enlargement of the category of Deligne $1$-motives  where duals may be introduced for commutative group schemes even in the above bad cases, and then introduces a new kind of topology on their cohomology groups which is necessary for introducing modifications in order to make dualities perfect. This is the bulk of the technical work; once the stage has been set, the proofs proceed by d\'evissage to known cases (as has been the case before for most results in the area). It should be noted that the dualities are not symmetric here like in many previous duality theorems.

There is no question that these results are nontrivial and technically innovative. Where I have some doubts is regarding their potential applicability. After all, arithmetic duality theorems have been invented in order to synthesize fundamental laws in arithmetic and to edict a framework that makes  them applicable to a number of problems. Here I don't know, for instance, how one could exploit the fact that $H^1$ of a commutative group scheme is dual to the maximal Hausdorff quotient of some complicated topological group.  For this reason I am not sure the paper is of interest to a wide circle of arithmetic geometers, but it may be appreciated by some experts."

\bigskip

I will now try to respond to the above comments.

\medskip

\section*{Notation}

\begin{spacing}{1.15}

If $k$ is a field, {\it ``algebraic over $k$"} means {\it of finite type over $k$}. The symbol $\ksep$ denotes a fixed separable algebraic closure of $k$.

A {\it $k$-variety} is a {\it separated, geometrically integral and algebraic $k$-scheme}.

If $X$ is a $k$-variety,
\[
\br_{\lbe 1}X\defeq\krn[\e\br X\to\br X_{\ksep}\e]
\]
and
\[
\hskip .7cm\br_{\be a}X\defeq\br_{\lbe 1}X/\le\img\be[\le\br k\to \br X\le].
\]
If $k$ is a global field,
\begin{equation}\label{bra}
\hskip 1.55cm \Bha(X)\defeq\krn\be[\e\br_{\be a}X\to\textstyle{\prod}_{\e{\rm all\,} v}\e\br_{\be a}X_{k_{v}}].
\end{equation}

If $B$ is an abelian group, $B^{\le *}\defeq \Hom\le(B,\Q/\Z\le)$.

If $G$ is an algebraic $k$-group, $Z_{G}$ denotes its (scheme-theoretic) center \cite[${\rm VI}_{\rm B}$, Corollary 6.2.5(i)]{sga3} and, if $G$ is smooth, $\s D(G\e)$ denotes its derived $k$-subgroup \cite[Definition A.1.14]{cgp}.

Below, we write $p$ for the characteristic of $k$ if $\car k>0$ and $H^{\le i}(k,-)$ denotes fppf (hyper)cohomology.

\section{Algebraic groups over global function fields.}\label{aggf}

As is well-known, the arithmetic theory of linear algebraic groups over number fields is a well-developed subject. By contrast, the corresponding theory over global function fields (especially in the non-reductive case) is only incipient, in spite of \cite{cf,ros18,ros20,dh22}. In this section I discuss certain open problems in the latter area and explain how the topological methods developed in \cite{ga23} are relevant to the eventual solution of these problems.

One of these problems involves the Hasse principle for torsors (or principal homogeneous spaces) under linear algebraic groups. Regarding this problem, the following theorem of Sansuc and Chernousov (see \cite[Theorem 8.5 and Corollary 8.7]{san} and \cite{che}) is well-known:

\begin{theorem}{\rm (Sansuc-Chernousov)}\label{sat} If $k$ is a number field and $G$ is a connected linear algebraic $k$-group, then the Tate-Shafarevich set $\Sha^{\le 1}\lbe(G\le)$ is equipped with a canonical abelian group structure and there exists a canonical perfect pairing of finite abelian groups
\begin{equation}\label{sp}
\Sha^{\le 1}\lbe(G\le)\times \Bha(G\le)\to\Q/\Z.
\end{equation}
Consequently, if $X$ is a torsor under $G$ over $k$, then the Brauer-Manin obstruction attached to $\Bha(X)$ is the only obstruction to the Hasse principle for $X$.
\end{theorem}

See \cite[\S5.2]{sko} for the basics on the Brauer-Manin obstruction.

\begin{remark}\label{abe} If $k$ is any global field and $G$ is a {\it smooth} connected linear algebraic $k$-group, then there exists a canonical (i.e., functorial) exact sequence of abelian groups (cf. \cite[Lemma 6.9(i)]{san}):
\[
\pic(G_{\lbe\ksep}\be)^{\lbe\g}\to H^{2}(\g,G^{\lle D}\be(\ksep))\to \br_{\lbe a}\le G\to H^{\le 1}\lbe(\g,\pic\e G_{\lbe\lle\ksep}\be)\to H^{3}(\g,G^{\lle D}\be(\ksep)),
\]	
where $\g={\rm Gal}\le(\ksep\be/k\lle)$ and $G^{\lle D}$ is the Cartier dual of $G$. Consequently, the second and third maps above induce maps of abelian groups
\begin{equation}\label{shm}
\Sha^{\e 2}\lbe(\g,G^{\lle D}\be(\ksep))\to \Bha(G\le)\to \Sha^{\e 1}\lbe(\g,\pic\e G_{\lbe\lle\ksep}\be).
\end{equation}
In certain cases, either the first or the second map above is an isomorphism. For example, if $\pic\e G_{\lbe\lle\ksep}=0$ (respectively, $G^{\lle D}\be(\ksep)=1$\,\footnote{\e This is the case, for example, if $G$ is unipotent since such groups have no nontrivial characters over a field \cite[XVII, Proposition 2.4(ii)]{sga3}.}), then the first (respectively, second) map in \eqref{shm} is an isomorphism. See Examples \ref{ex1} and \ref{ex2} below.
\end{remark}

The proof of Theorem \ref{sat} is easily reduced to the reductive case and depends crucially (via the additivity lemma \cite[Lemma 6.6]{san}) on the $\ksep$-rationality of reductive algebraic $k$-groups. Now, if $k$ is a global {\it function} field, then Theorem \ref{sat} remains valid if $G$ is a {\it reductive} algebraic $k$-group\,\e\footnote{ By \cite[Corollary 5.10]{gaq}, \cite[Theorem 5.2]{dh20} and \cite[Theorem 3.7]{dh22}. I also note that the last assertion in Theorem \ref{sat} for reductive groups over global function fields is now a particular case of \cite[Theorem 2.5]{dh22}.}\,, but fails if $G$ is an arbitrary non-reductive smooth and connected linear algebraic $k$-group. Non-reductive groups, even commutative ones, are much more complicated than reductive groups. For example, there exist commutative (non-reductive) pseudo-reductive  \cite[Definition 1.1.1]{cgp} algebraic $k$-groups which are not even $\ksep$-unirational \cite[Example 11.3.1]{cgp}. Nevertheless, there exist interesting {\it specific classes} of non-reductive linear algebraic $k$-groups whose arithmetic theory should, in a certain sense, be ``sufficiently close" to that of reductive $k$-groups to justify expectations of rapid progress in this area for such groups.

\subsection{Groups of Lourdeaux type} \label{glt} One of these\,\e\footnote{\e For other classes, see Remark \ref{rmn}(b) and subsection \ref{qrg} below.} is the class of {\it groups of Lourdeaux type}, i.e., algebraic $k$-groups of the form $\rkk(G^{\le\prime}\le)/\rkk(Z^{\le\prime}\le)$, where $\kp\be/k$ is a nontrivial finite and purely inseparable field extension, $G^{\le\prime}$ is a 
reductive algebraic $\kp$-group and $Z^{\le\prime}$ is a (possibly non-smooth) central $\kp$-subgroup of $G^{\e\prime}$. These pseudo-reductive $k$-groups were studied in \cite{lord} (over any field $k$) and shown there to be {\it retract $\ksep$-rational}, which implies that they {\it satisfy Sansuc's additivity lemma \cite[Lemma 6.6]{san}} (see \cite[Proposition 2.4]{lord}). Consequently, a substantial part of the theory developed for reductive groups in \cite{san} (and in several other works where \cite[Lemma 6.6]{san} plays a central role, e.g., \cite{bvh}) carries over, {\it mutatis mutandis}, to the class of groups of Lourdeaux type. For example, if $G$ is a group of Lourdeaux type, then the pairing \ref{sp} exists (independently of any abelian group structure on $\Sha^{\le 1}\lbe(G\le)$). See \cite[Proposition 8.1]{san}. Now two basic questions arise: if $G$ is a group of Lourdeaux type over a global function field $k$, does $\Sha^{\le 1}\lbe(G\le)$ admit an abelian group structure (canonical or otherwise)? And, if so, is the pairing \eqref{sp} perfect? I.e., is the Brauer-Manin obstruction attached to $\Bha(X)$ the only obstruction to the Hasse principle for torsors $X$ under $G\e$? In the case of reductive $k$-groups, the known (affirmative) answers to these questions (which are now particular cases of results in \cite{dh22}) make crucial use of certain ``abelianization techniques", whose basic ingredient is Deligne's map
\begin{equation}\label{dmap}
\rho\colon \widetilde{G}\onto \s D(G\e)\into G,
\end{equation}
where  $\widetilde{G}$ is the simply-connected central cover of $\s D(G\e)$. For non-reductive $k$-groups $G$, the map \eqref{dmap} is unavailable. However, Conrad and Prasad have introduced an analog of the simply-connected central cover of a connected semisimple $k$-group \cite[Chapter 5]{cp} which suffices for developing useful abelianization techniques for a broad class of (not necessarily affine) smooth and connected algebraic $k$-groups. These abelianization techniques (which are reviewed below) apply to all pseudo-reductive $k$-groups, hence in particular to all $k$-groups of Lourdeaux type.

\subsection{Abelianization techniques for a class of smooth and connected algebraic groups}\label{abt} Let $k$ be any field and let $G$ be a smooth and connected algebraic $k$-group. By \cite[${\rm VI}_{\le\rm B}$, Propositions 7.1(i) and 7.8]{sga3}, the derived $k$-subgroup $\s D(G\le)$ of $G$ is smooth and connected. Consider the exact and commutative diagram of algebraic $k$-groups 
\begin{equation}\label{diag}
\xymatrix{1\ar[r]& L\cap G_{\rm sab}\ar@{^{(}->}[d]\ar[r]&G_{\rm sab}\ar@{^{(}->}[d]\ar[r]& A\ar@{=}[d]\ar[r]&1\\
1\ar[r]& L\ar@{->>}[d]\ar[r]&G\ar@{->>}[d]\ar[r]&A\ar[r]&1\\
&G^{\rm aff}\ar@{=}[r]&G^{\rm aff}&&
}
\end{equation}
where the middle row (respectively, column) is the Albanese (respectively, anti-Chevalley) decomposition of $G$. Clearly $\s D(G\le)\subseteq L=\krn\alb_{\le G}$, whence $\s D(G\le)$ is affine. Now, since $G_{\rm sab}$ is central in $G$ and $G=G_{\rm sab}\cdot L$ \cite[Theorem 5.1.1(1)]{bri17a}, we have $Z_{L}=Z_{G}\cap L$ and $\s D(G\le)=\s D(L)$. In particular, $\s D(L)$ is smooth (even if $L$ is not) and $Z_{G}\cap \s D(G\le)=Z_{L}\cap \s D(L)$. Further, diagram \eqref{diag} induces an exact and commutative diagram of commutative algebraic $k$-groups (cf. \cite[(1.9)]{ga23})
\begin{equation}\label{diag2}
\xymatrix{0\ar[r]& L\cap G_{\rm sab}\ar@{^{(}->}[d]\ar[r]&G_{\rm sab}\ar@{^{(}->}[d]\ar[r]& A\ar@{=}[d]\ar[r]&0\\
0\ar[r]& Z_{L}\ar[d]\ar[r]&Z_{G}\ar[d]\ar[r]&A\ar[r]&0\\
&Z_{G^{\rm aff}}\ar@{=}[r]&Z_{G^{\rm aff}}.&&
}
\end{equation}

\medskip

Now assume that

\smallskip

\begin{enumerate}
\item[\textbf{(H)}] $\s D(L)$ is {\it perfect}, i.e., $\s D(\s D(L\le))=\s D(L\le)$.
\end{enumerate}

\medskip

By \cite[Proposition 1.2.6]{cgp}, \textbf{(H)} holds if $L$ is a pseudo\le-\lle reductive $k$-group.

\smallskip

Since $\s D(L)$ is smooth, connected, affine and perfect, we may consider its universal smooth $k$-tame central cover $\widetilde{\s D(L)}$ \cite[\S5.2]{cp}, which will henceforth be denoted by $\widetilde{L}$. Thus we have a central extension of affine algebraic $k$-groups
\begin{equation}\label{cc}
1\to Z\to \widetilde{L}\overset{\!q}{\to} \s D(L)\to 1,
\end{equation}
where $Z$ is {\it $k$-tame}, i.e., contains no non-trivial unipotent $k$-subgroup schemes. Now set
\begin{equation}\label{gz}
G_{\lle 0}=Z_{G}\s D(L)=Z_{G}\s D(G\lle),
\end{equation}
which is a closed normal $k$-subgroup of $G$ such that $G/G_{0}$ is commutative. Now let 
\begin{equation}\label{pdm}
\rho\colon \widetilde{L}\overset{\!q}{\onto} \s D(L)\into G_{\lle 0}
\end{equation}
be Deligne's map (note that, if $G=L$ is reductive, then 
$G_{\lle 0}=Z_{L}\s D(L)=L$ by \cite[XXII, 6.2.3]{sga3} and \eqref{pdm} agrees with \eqref{dmap}). We let $\s D(L)$ (respectively, $Z_{G}$) act on $\widetilde{L}$ via  $q$ ``\le by conjugation" (respectively, trivially). Since the map $Z_{\widetilde{L}}\to Z_{\s D(L)}$ induced by $q$ \eqref{cc} is surjective \cite[Proposition 2.2.12(2)]{cgp}, the $2\e$-term complex\,\e\footnote{ Of fppf sheaves of groups on $k$ in degrees $-1$ and $0$.} $[\e\widetilde{L}\overset{\!\rho}{\to}G_{\lle 0}\e]$ is a {\it quasi-abelian crossed module} \cite[Definition 3.2]{gaq}. Thus, by \cite[Proposition 3.4]{gaq}, $[\e\widetilde{L}\overset{\!\rho}{\to} G_{\lle 0}\e]$ is quasi-isomorphic to the complex of commutative algebraic $k$-groups
\begin{equation}\label{cg}
\hskip 1.8cm C_{G_{\lle 0}}\defeq [\e Z_{\lle\widetilde{L}}\to Z_{G_{0}}\e]= [\e Z_{\lle\widetilde{L}}\to Z_{G}\le Z_{\s D(\lbe L)}\e].
\end{equation}
For every integer $i\geq -1$, set
\begin{equation}\label{abc}
H^{\le i}_{\rm ab}(G_{0})\defeq H^{i}(k,C_{G_{\lle 0}}\lbe).
\end{equation}
By \cite[sequences (3.4) and (3.5)]{gaq}, the abelian group $H^{\le i}_{\rm ab}(G_{\lle 0}\e)$ fits into canonical exact sequences of abelian groups:  
\[
\hskip 1.5cm\dots\to H^{\e i}(Z_{G}Z_{\s D(L)})\to H^{\le i}_{\rm ab}(G_{\lle 0}\e)\to H^{\e
i+1}(Z_{\widetilde{L}})\to H^{\e i+1}(Z_{G}Z_{\s D(L)})\to\dots.
\]
and 
\[
\hskip .4cm\dots\to H^{\e i-1}(Z_{G}/Z_{L}\le\cap\le \s D(L))\to H^{\e
i+1}(Z\lle)\to H^{\le i}_{\rm ab}(G_{\lle 0}\e)\to H^{\e i}(Z_{G}/Z_{L}\le\cap\le \s D(L))\to\dots, 
\]
where $H^{\e i}(-)\defeq H^{\e i}(k,-)$ and $Z=\krn q$ \eqref{cc}. Further, by \cite[Theorem 4.2]{gaq}, there exists a commutative diagram of pointed sets with exact row and column
\begin{equation}\label{dou}
\hskip 1.5cm\xymatrix{&H^{\e 1}(\widetilde{L}\e)\ar[dr]\ar[d]&&\\	
H^{\e 0}(G/G_{\lle 0})\ar[dr]\ar[r]&H^{\e 1}(G_{\lle 0})\ar[d]^(.45){{\rm ab}^{\lbe 1}}\ar[r]&H^{\le 1}\lbe(G\le)\ar[r]^(.45){\pi^{(1)}}&H^{\e 1}(G/G_{\lle 0})\\
&H^{\le 1}_{\rm ab}(G_{\lle 0})\ar[d]&&\\
&H^{\e 2}(\widetilde{L}\e),
}
\end{equation}
where $\pi^{(1)}$ is induced by the projection $\pi\colon G\onto G/G_{\lle 0}$ and the column is induced by the exact sequence of complexes
\[
1\to [\e 1\to G_{\lle 0}\e]\to [\e\widetilde{L}\to G_{\lle 0}\e]\to [\e\widetilde{L}\to 1\e]\to 1
\]
via the quasi-isomorphism $[\e\widetilde{L}\to G_{\lle 0}\e]\simeq C_{\lbe G_{\lle 0}}$\,\e\footnote{\e In diagram \eqref{dou}, exactness of the column at $H^{\le 1}_{\rm ab}(G_{\lle 0})$ is meant in the sense of Douai, i.e., $y\in H^{\e 1}_{{\rm{ab}}} (G_{\lle 0}\e)$ is in the image of ${\rm{ab}}^{1}$ if, and only if, the image of $y$ in (Giraud's) $H^{\e 2}(\widetilde{L}\e)$ is {\it neutral}.}\,. Note that the pointed sets $H^{\le i}(G/G_{\lle 0})$ in \eqref{dou} are, in fact, abelian groups. 

\begin{remarks}\indent
\begin{enumerate}
\item[(a)] If $L$ is {\it reductive}, so that $L=Z_{L}\s D(L)$, then
\[
G=G_{\rm sab}\e L\subseteq Z_{G}Z_{L}\s D(L)=Z_{G}\s D(L)=G_{0}\subseteq G,
\]
i.e., $G_{\lle 0}=G$. Thus the map $\pi^{(1)}$ in diagram \eqref{dou} is the zero map and the complex \eqref{cg} is $C_{G_{\lle 0}}=[\e Z_{\lle\widetilde{L}}\to Z_{G}\e]$. Consequently,  the above setting generalizes that which is described in \cite[beginning of \S4.1, pp.~113\e-\lle 115]{d11}.

\item[(b)] Condition \textbf{(H)} certainly fails for non\le-\le pseudo\le-\le reductive algebraic $k$-groups, such as the non\le-\le commutative $k$-wound group described in \cite[Example B.2.9]{cgp}. Thus, the abelianization techniques described above are insufficient for studying the arithmetic of general affine and connected algebraic $k$-groups. Nevertheless, certain specific classes of non\le-\le pseudo\le-\le reductive $k$-groups can be studied via other methods. See subsection \ref{qrg} below for a brief discussion of certain arithmetic aspects of an interesting class of (non\le-\le pseudo\le-\le reductive) quasi\le-\le reductive algebraic $k$-groups.
\end{enumerate}
\end{remarks}

\smallskip

We now let $k$ denote a global field.

\smallskip

\subsubsection{Affine groups} Assume that $A=0$ in diagram \eqref{diag}, i.e., $G=L$ is smooth, affine and connected.

If $G$ is {\it reductive} (so that, as noted above, $G_{0}=Z_{G}\s D(\lbe G\e)=G\e$), the complex 
$C:=[\e Z_{\lle\widetilde{G}}\to Z_{\lle G}\e]$ \eqref{cg} plays a central role in the arithmetic theory of $G$ via isomorphisms $\Sha^{\le 1}\lbe(G\le)\simeq \Sha^{\le 1}\lbe(C\le)$ and $\br_{\be a}\le G\simeq H^{\le 1}(k,C^{\lle D})$\,\footnote{\e See \cite[Theorem 5.13]{bor}, \cite[Corollary 5.10]{gaq} and \cite[Theorem 3.7]{dh22}.}\,. A key result in this setting is the global duality theorem $\Sha^{\le 1}\lbe(C\le)\simeq \Sha^{\le 1}\lbe(C^{D}\le)^{*}$ \cite[Theorem 5.2]{dh20}, which is equivalent to the duality statement in Theorem \ref{sat} for a reductive group over an arbitrary global field. The proof of the aforementioned duality theorem makes use of the fact that the complex $C$ above is quasi-isomorphic to a $2$-term complex of tori\,\e\footnote{ This is, essentially, a consequence of the fact that, in a reductive group, Cartan subgroups and maximal tori coincide \cite[\S13.17, Corollary 2(c), p.~175]{bo}.}\,, so that the theorem in question is, in fact, a global duality theorem for $2$-term complexes of $k$-tori. Now, in order to obtain such a global duality theorem, one first needs to establish local (Pontryagin) duality theorems involving the canonical pairings of abstract abelian groups
\begin{equation}\label{lpair}
H^{\le i}(k,C^{\lle D}\le)\times H^{1-i}(k,C\e)\to \br k\simeq \Q/\Z,
\end{equation}
where $k$ is a non-archimedean local field and $i=0,1$ or $2$. A key requirement to achieve the latter is to equip the abelian groups $H^{\le i}(k,C^{D}\le)$ and $H^{1-i}(k,C\le)$ with appropriate locally compact topologies. When $C=[\e T_{1}\to T_{2}\e]$ is a $2$-term complex of $k$-tori, this goes as follows (cf. \cite[beginning of \S3]{dem11}): consider the exact sequence of abstract abelian groups
\begin{equation}\label{tseq}
H^{\lle i}(k,T_{1}^{\lle D})\overset{\!f^{(i)}}{\lra} H^{\lle i}(k, C^{\lle D}\le)\overset{\!g^{\lle(i)}}{\lra}  H^{\le i+1}(k,T_{2}^{\lle D}).
\end{equation}
If $i=0$ or $1$, the natural topology to place on $H^{\lle i+1}(k,T_{2}^{D})$ is the {\it discrete topology} (in fact, $H^{\lle 1}(k,T_{2}^{D})$ is finite and discrete). Consequently, the Nagao topology on $H^{\lle i}(k, C^{D}\le)$ determined by the sequence \eqref{tseq} and a choice of a set-theoretic section of the map $g^{\lle(i)}$ in \eqref{tseq} is {\it independent of this choice} and, therefore, {\it uniquely determined} \cite[Remark 2.4]{ga23}. Since $H^{\lle 1}(k,T_{1})$ is also (finite and) discrete, $H^{\le i}(k,C\le)$ can be similarly equipped with a uniquely determined Nagao topology. Now \eqref{lpair} is a pairing of locally compact abelian topological groups which induces the perfect pairings stated in \cite[Theorem 3.1]{dem11}. The local perfect pairings alluded to above are basic ingredients in the proof of the global duality statement $\Sha^{\le 1}\lbe(C\le)\simeq \Sha^{\le 1}\lbe(C^{D}\le)^{*}$, which itself is one of the ingredients needed to construct the Poitou-Tate exact sequences for $2$-term complexes of tori established in \cite{dh20}. These Poitou-Tate exact sequences, namely \cite[Theorems 5.7 and 5.10]{dh20}, play a central role in the study of weak and strong approximation for homogeneous spaces under reductive groups with {\it reductive stabilizers} over a global function field \cite[Theorems 4.2 and 5.8]{dh22}.

Now, if $k$ is a global function field and $G$ is a {\it non\le-\le reductive} smooth, affine and connected algebraic $k$-group which satisfies condition \textbf{(H)}, i.e., such that $\s D(\lbe G\e)$ is perfect, then the complex $C=C_{G_{0}}=[\e Z_{\lle\widetilde{G}}\to Z_{G}Z_{\s D(G\le)}\le]$ \eqref{cg} is not, in general, quasi-isomorphic to a $2$-term complex of $k$-tori. Thus, in order to extend the main results of \cite{dh22} to some interesting classes of non-reductive algebraic $k$-groups (e.g., groups of Lourdeaux type) using similar techniques, one needs to establish Poitou\e-Tate exact sequences for $2$-term complexes of {\it arbitrary} affine and commutative algebraic $k$-groups\,\e\footnote{\e One also has to contend with the fact that, in this non-reductive setting, there exists an additional abelianization map, namely the map $\pi^{(1)}$ in diagram \eqref{dou}. For lack of time, I cannot discuss this problem here.}\,. To this end, the first task is to establish local duality theorems involving the pairings \eqref{lpair} for such complexes $C$. As before, the latter requires equipping the abelian groups $H^{\le i}(k,C^{\lle D}\le)$ and $H^{1-i}(k,C\le)$ with (possibly non-Hausdorff) locally compact topologies, where $i=0,1$ or $2$, $k$ is a local function field and $C=[\e C_{1}\to C_{2}\le]$ is a $2$-term complex of affine and commutative algebraic $k$-groups. In the case $C_{1}=0$ above, i.e., $C=C_{2}$ is a commutative affine algebraic $k$-group concentrated in degree $0$, Rosengarten \cite[\S3]{ros18} equipped the groups $H^{\le i}(k,C_{2}^{D}\le)$ and $H^{1-i}(k,C_{2}\le)$ with suitable Hausdorff and locally compact topologies and established Pontryagin duality theorems for these groups (see \cite[Theorem 4.2]{ga23} for a summary of the corresponding statements). 
Now, in order to topologize the groups $H^{\le i}(k,C^{\lle D}\le)$ and $H^{1-i}(k,C\le)$ for a general complex $C=[\e C_{1}\to C_{2}\le]$ as above, we proceed as before, i.e., we consider the exact sequence of abstract abelian groups
\begin{equation}\label{cseq}
H^{\lle i}(k,C_{1}^{\lle D})\overset{\!f^{(i)}}{\lra} H^{\lle i}(k, C^{\lle D}\le)\overset{\!g^{\lle(i)}}{\lra}  H^{\le i+1}(k,C_{2}^{\lle D})
\end{equation}
and equip $H^{\lle i}(k, C^{\lle D}\le)$ with any one of the Nagao topologies induced by the above sequence and a choice of a set-theoretic section $\sigma_{i}$ of the map $g^{\lle(i\lle)}$ in \eqref{cseq}, where the groups $H^{\lle i}(k,C_{1}^{\lle D})$ and $H^{\lle i+1}(k,C_{2}^{\lle D})$ are equipped with their respective Rosengarten topologies. The topology of $H^{1-i}(k,C\le)$ is defined similarly and depends on the choice of a set-theoretic section $\tau_{i}$ of the map $H^{\lle 1-i}(k, C\le)\to H^{\le 2-i}(k,C_{1})$ induced by the canonical morphism of complexes $C\to C_{1}[1]$. By the work done in \cite{ga23}, it is reasonable to expect that the topologies of the associated maximal Hausdorff quotients $H^{\lle i}_{\be\sigma_{\lbe i}}\be(k, C^{\lle D}\le)_{\lle\rm Haus}$ and $H^{\le 1-i}_{\tau_{\lbe i}}(k, C\le)_{\lle\rm Haus}$ will be independent of the choices of $\sigma_{i}$ and $\tau_{i}$. Further, these groups should satisfy (appropriate) Pontryagin duality theorems.

\begin{remark}\label{dev} I now stress an important technical point. In establishing the preceding (interdependent) assertions, namely, independence of the topologies on the choice of set-theoretic sections {\it and} the corresponding local duality theorems, the  ``functoriality lemma" \cite[Lemma 2.1]{ga23} {\it must} play a central role. The reason is that, in a setting where the use of d\'evissage techniques is necessary and the topologies to be placed on various abstract abelian groups {\it depend} on the choice of d\'evissage(s), a statement such as \cite[Lemma 2.1]{ga23} is {\it essential\e} for verifying the continuity of the relevant pairings (cf. \cite[proof of Proposition 4.10]{ga23}). Consequently, \cite[Lemma 2.1]{ga23} should (in particular) play a key role in extending to general $k$-$1$-motives the local duality theorems for pure $k$-$1$-motives established in \cite{ga23} (using the d\'evissage that yields \cite[Theorem 2.3]{hsz1}). Unfortunately, the topological verifications alluded to above have been overlooked in several papers that claim to have proven a local duality theorem. This is the case, for example, with \cite{vhd}, where the indicated verifications have been omitted throughout, e.g., in the proofs of \cite[Proposition 4.1 and Theorem 4.3]{vhd}. Very regrettably, the author of \cite{vhd} is no longer with us to address this issue himself. In his stead, Mr. Rivera-Mesa (a student of Lucchini Arteche's whom I'm co\le-\le advising) will include all the required topological details in his forthcoming generalization of \cite{vhd} to certain classes of singular varieties over $p\e$-adic fields.
\end{remark}

Now assume, for the sake of argument, that the appropriate local duality theorems for $2$-term complexes $C=[\e C_{1}\to C_{2}\le]$ as above have been established. 
These local duality theorems should then be the basic ingredients that lead to the construction of suitable {\it Poitou\le-Tate exact sequences} for complexes $C$ as above over global function fields. Such
Poitou\le-Tate exact sequences have already been established in two particular cases, namely by Demarche and Harari when $C_{1}$ and $C_{2}$ are $k$-tori \cite[Theorems 5.7, 5.8 and 5.10]{dh20} and by Rosengarten when $C_{1}=0$ and $C_{2}$ is arbitrary \cite[Theorems 5.17.1 and 5.17.2]{ros18}. In each of these particular cases, the indicated Poitou-Tate exact sequences have had important consequences for the arithmetic theory of reductive algebraic $k$-groups \cite[Theorems 2.5 and 4.2]{dh22} and connected algebraic $k$-groups that are either commutative or pseudo-reductive \cite[Theorems 1.1 and 1.6]{ros20}. Thus, it is reasonable to expect that Poitou-Tate exact sequences for {\it arbitrary} $2$-term complexes of affine and commutative algebraic $k$-groups {\it will} have important consequences for the arithmetic theory of smooth, affine and connected algebraic $k$-groups whose derived subgroup is perfect.

\subsubsection{Non-affine groups} The arithmetic theory of arbitrary (i.e., not necessarily linear) algebraic groups has been advanced by a number of authors, but (to my knowledge) only in the number field case. One of the problems discussed in the literature is the following question:

\begin{question}\label{q1}\cite[p.~133]{sko}\, If $G$ is a connected algebraic group over a number field $k$ and $X$ is a torsor under $G$ over $k$, is the Brauer-Manin obstruction attached to $\Bha(X)$ \eqref{bra} the only obstruction to the Hasse principle for $X$?
\end{question}

By Theorem \ref{sat}\,, the above question has an affirmative answer if $G$ is linear. However, the answer is negative in general. Indeed, in \cite[\S3.6]{bcts} a non-linear and non\e-\e commutative $k$-group $G$ is constructed for which the answer to Question \ref{q1} is negative. But what about {\it commutative} (non-linear) algebraic groups $G\e$ over a number field $k$? Since $\car k=0$, we may restrict our attention to semiabelian $k$-varieties $G$, i.e., extensions 
\begin{equation}\label{sab}
0\to T\to G\to A\to 0,
\end{equation}
where $T$ is a $k$-torus and $A$ is an abelian $k$-variety. This reduction step follows from the fact that, over a field $k$ of characteristic $0$, every commutative unipotent algebraic $k$-group is a $k$-vector group (i.e., is isomorphic to $\G_{\lbe a,\le k}^{n}$ for some $n\geq 0$) \cite[IV,\S2, 4.2, p.~498]{dg}. Now, since the answer to Question \ref{q1} is affirmative for $T$ by Theorem \ref{sat} and also for $A$ if $\Sha^{1}(A)$ is finite by
\cite[Theorem 6.2.3, p.~125]{sko}, it is reasonable to expect that Question \ref{q1} 
also has an affirmative answer for semiabelian $k$-varieties $G$ \eqref{sab} (under the stated finiteness hypothesis). This is, indeed, the case by work of Harari and Szamuely \cite[Theorem 1.1]{hsz2}. However, the proof of \cite[Theorem 1.1]{hsz2} is not a straightforward generalization of the previously known proofs in the particular cases $G=T$ and $G=A$ alluded to above. In fact, the proof doesn't even take place inside the category of semiabelian $k$-varieties. Duality considerations force the authors of \cite{hsz2} to work in the (larger) category $\s{M}_{k,\e 1}^{\le\rm Del}$ of Deligne $k$-$1$-motives. More precisely, the proof of \cite[Theorem 1.1]{hsz2} depends crucially on the global duality theorem for such objects previously established in \cite[Theorem 0.2]{hsz1} (applied to $G$ \eqref{sab} and its dual Deligne $k$-$1$-motive $G^{\le\vee}=[\e T^{\lle D}\!\!\to\!\lbe A^{\lle t}\e]\e$).

\smallskip

What about the analog of the above problem over {\it global function fields}\le? 
As noted above, Question \ref{q1} has an affirmative answer for abelian varieties with finite Tate-Shafarevich group by \cite[Theorem 6.2.3, p.~125]{sko}. Since the indicated result is also valid in the function field case (with essentially the same proof), a natural version of (the commutative case of) Question \ref{q1} over a global function field is

\begin{question} \label{q2} Let $k$ be a global function field and let $G$ be a smooth, connected and commutative algebraic $k$-group with Albanese decomposition
\begin{equation}\label{lgb}
0\to L\to G\to A\to 0,
\end{equation}
where $L$ is affine and connected, $A$ is an abelian $k$-variety and $\Sha^{1}\lbe(\be A\lle)$ is finite. Assume that, for every torsor $X$ under $L$ over $k$, the Brauer-Manin obstruction attached to $\Bha(X)$ \eqref{bra} is the only obstruction to the Hasse principle for $X$, i.e., $X(k)=\emptyset\iff X(\A_{k})^{\lbe\Bha(X)}=0$. Is the corresponding statement with $G$ in place of $L$ also true?
\end{question}

The above question is {\it much more difficult} than the commutative case of Question \ref{q1}\,, chiefly because in Question \ref{q2} we can no longer restrict our attention to semiabelian $k$-varieties $G$. Indeed, and in stark contrast with the number field case, the class of commutative unipotent algebraic groups over a global function field is much larger (and much more complicated) than the class of $k$-vector groups. Nevertheless, we can still try to answer Question \ref{q2} via a suitable generalization of the Harari-Szamuely method developed in \cite{hsz2}. To do that, a first step is to enlarge the category $\s{M}_{k,\e 1}^{\le\rm Del}$ of Deligne $k$-$1$-motives and work in the new category $\s{M}_{k,\e 1}$ introduced in \cite{ga23} (which, by construction, contains the $k$-$1$-motive $G^{\le\vee}=[\e L^{\be D}\!\!\to\!\lbe A^{\lle t}\,]\e$ dual to $G$ \eqref{lgb}). Further steps include proving appropriate duality theorems in $\s{M}_{k,\e 1}$ that generalize those obtained in \cite{hsz1} for Deligne $1$-motives over number fields. 

\smallskip

My interest in Question \ref{q2} was one of the reasons for introducing the new theory of $k$-$1$-motives described in \cite{ga23}.

\smallskip

\begin{example}\label{ex1}
Let $\kp\!/k$ be a nontrivial finite and purely inseparable field extension and let $G$ be a semiabelian $k$-variety which is not an abelian variety, with Albanese decomposition $0\to T\to G\to A\to 0$ (thus $T\neq 0$). By \cite[Corollary A.5.4(3), p.~508]{cgp}, the above sequence induces an extension of smooth, connected and commutative algebraic $k$-groups
\begin{equation}\label{sub}
0\to \rkk(\le T_{\kp}\be)\to \rkk(G_{\kp}\be)\to  \rkk(A_{\kp}\be)\to 0.
\end{equation}
The pullback of \eqref{sub} along the embedding $A\into \rkk(A_{\kp}\be)$ is an extension of the form \eqref{lgb}, namely
\begin{equation}\label{lgb2}
0\to \rkk(T_{\kp}\le)\to H\to A\to 0,
\end{equation}
where $H=\rkk(G_{\kp}\le)\times_{\rkk(A_{\kp}\lle)}A$. Set $L=\rkk(\le T_{\kp}\be)$.
Since $L$ is $\ksep$-rational, $L$ satisfies \cite[Lemma 6.6]{san}, whence the pairing $\Sha^{1}(L)\times \Bha(L\le)\to \Q/\Z$ \eqref{sp} exists. We have $\Sha^{1}(L)=\Sha^{1}(\kp, T\e)=\Sha^{1}(T\e)$ (by the smoothness of $T$ and the topological invariance of the \'etale site \cite[O4DY, Proposition 59.45.4(3)]{sp}). On the other hand, $\pic\e L_{\lbe\lle\ksep}=0$ by \cite[\S4.1]{ach19a}, whence the map $\Sha^{2}(\g,L^{\be D}\be(\ksep))\to \Bha(L\le)$ in \eqref{shm} is an isomorphism. Now, since $k$ has $p\e$-cohomological dimension $1$, \cite[Theorem A.4]{lord} yields the equality $\Sha^{2}\lbe(\g,L^{\be D}\be(\ksep))=\Sha^{2}(\g,T^{\lle D}\be(\ksep))=\Sha^{2}(T^{\lle D})$. Thus the pairing $\Sha^{1}(L)\times \Bha(L\le)\to \Q/\Z$ \eqref{sp} is isomorphic to a pairing $\Sha^{1}(T\e)\times \Sha^{\lle 2}\lbe(T^{\lle D})\to \Q/\Z$ which should coincide (up to a sign) with the (perfect) Tate\e-\le Nakayama pairing $\Sha^{1}(T\e)\times \Sha^{\lle 2}\lbe(T^{\lle D})\to \Q/\Z$. Thus, up to certain verifications, the commutative pseudo\e-\le reductive $k$-group $L=\rkk(T_{\kp})$ satisfies the condition stated in Question \ref{q2}\,, i.e., if $X$ is a torsor under $L$ over $k$, then the Brauer\e-\le Manin obstruction attached to $\Bha(X)$ is the only obstruction to the Hasse principle for $X$. Question \ref{q2} asks whether or not a similar statement is valid for the extension $H$ \eqref{lgb2}.
\end{example}

\begin{example} \label{ex2} Let the setting be that of the previous example. Then $T$ is the maximal $k$-torus in $L=\rkk(\le T_{\kp}\be)$. Set $U=\rkk(\le T_{\kp}\be)/\lle T$ and let $\overline{H}=H/\lle T$. The sequence \eqref{lgb2} induces an exact sequence of smooth, connected and commutative algebraic $k$-groups
\begin{equation}\label{lgb3}
0\to U\to \overline{H}\to A\to 0.
\end{equation}
The $k$-group $U$ is a $k$-unirational unipotent $k$-wound group \cite[VI, \S5.1, Lemma, p.~70]{oes}. Thus, by \cite[Theorem 2.8(iii)]{ach19b}, $U$ {\it satisfies Sansuc's additivity lemma} \cite[Lemma 6.6]{san}. Consequently, the pairing
$\Sha^{\le 1}\lbe(U\le)\times \Bha(U\le)\to\Q/\Z$ \eqref{sp} is defined. 
Now, by Remark \ref{abe}\,, the map $\Bha(U\le)\to \Sha^{\e 1}\lbe(\g,\pic\e U_{\lbe\lle\ksep}\be)$ in \eqref{shm} is an isomorphism. Thus \eqref{sp} is isomorphic to a pairing $\Sha^{\le 1}\lbe(U\le)\times \Sha^{\le 1}\be(\g,\pic\e U_{\lbe\lle\ksep}\be)\to\Q/\Z$. It seems likely that, by suitably combining the work of Oesterl\'e \cite[VI, \S5.1-5.2, pp.~70\e-71]{oes}, Achet \cite[\S4.1]{ach19a} and Tate\e-\le Nakayama duality for $T$, the preceding pairing can be shown to be perfect. If this is indeed the case, then we can ask Question \ref{q2} for the extension $\overline{H}$ \eqref{lgb3}.
\end{example}

\begin{example}\label{ex3} Let $\kp\!/k$ be as in Example \ref{ex1} and let $B\onto A$ be a nontrivial infinitesimal isogeny of abelian $k$-varieties with kernel $I$. Then $G=\rkk(B_{\kp}\be)\times_{\rkk(A_{\kp}\be)}A$ is a pseudo\e-\e abelian $k$-variety with Albanese decomposition 
\[
0\to \rkk(\lbe I_{\kp}\!)\to G\to A\to 0.
\]
The (positive\e-\e dimensional\le) $k$-group $L=\rkk(I_{\kp}\be)$ is {\it smoothless}, i.e., $L^{\lbe\sm}=0$. Thus $\Sha^{1}(L\le)=\Sha^{1}(L^{\lbe\sm})=0$ by \cite[Example C.4.3]{cgp}, whence $L$ satisfies (vacuously) the condition stated in Question \ref{q2}\e. Thus we may ask Question \ref{q2} for $G$. 
\end{example}

\begin{remarks}\indent
\begin{enumerate}
\item[(a)] As noted previously, by a function field analog of the proof of \cite[Theorem 6.2.3, p.~125]{sko}, Question \ref{q2} has an affirmative answer for abelian varieties with finite Tate-Shafarevich group. Since the indicated proof depends on Cassels\e-\le Tate duality, it is natural, for the purpose of answering Question \ref{q2} for pseudo\e-\le abelian varieties, to seek to extend Cassels\e-\le Tate duality to such varieties (using their Albanese decomposition). This was another motivation for introducing the category of $k$-$1$-motives described in \cite{ga23}. Regarding this point, I should note that pseudo\e-\le abelian varieties $G$ were introduced in \cite{to13} and studied there mainly in terms of their anti\e-\le Chevalley decomposition $0\to B\to G\to U\to 0$, where $B$ is the maximal abelian $k$-subvariety of $G$ and $U$ is the corresponding unipotent quotient \cite[Theorem 2.1]{to13}. Unfortunately, I don't know how to develop a duality theory for pseudo\e-\le abelian $k$-varieties $G$ (or any other type of commutative algebraic $k$-groups, for that matter) using their anti\e-\le Chevalley decompositions.
		
\item[(b)] Smoothless commutative\,\e\footnote{\e See \cite[pp.~652\e-\le 653]{cgp} for an  interesting example of a {\it non\e-\le commutative} smoothless algebraic $k\le$-group.} algebraic $k$-groups, such as the $k$-group $\rkk(I_{\kp}\be)$ in Example \ref{ex3}\,, are important objects. For instance, every $k$-group of Lourdeaux type $\rkk(G^{\le\prime}\le)/\rkk(Z^{\le\prime}\le)$ (see subsection \ref{glt}) can be written as $\rkk(G^{\le\prime\prime}\le)/\rkk(\e\mu^{\le\prime}\le)$, where $G^{\le\prime\prime}\defeq G^{\le\prime}\lbe/\lbe\le(Z^{\le\prime}\le)^{\rm sm}$ is reductive and $\mu^{\le\prime}\defeq Z^{\le\prime}\lbe/\lbe(Z^{\le\prime})^{\rm sm}\subset G^{\le\prime\prime}$ is infinitesimal and central (and therefore of multiplicative type), so that $\rkk(\e\mu^{\le\prime}\le)$ is smoothless and central in $\rkk(G^{\le\prime\prime}\le)$. Further, by the comments after \cite[Definition 5.1.1, p.~66]{cp}, extensions of the form \eqref{cc} can have a smoothless $k$-tame kernel. Also, smoothless commutative algebraic groups appear in the problem of resolution of singularities over an imperfect field as Hironaka group schemes of a certain type (see \cite[Example 2.1]{miz73} and \cite[\S3]{oda73}), and play a central role in determining the structure of the Picard scheme of a proper variety over such a field. Additionally, smoothless commutative algebraic groups (which are necessarily affine) are closely related to the {\it projective} and geometrically non\le-\le reduced $k$\le-\le schemes studied by Scully \cite{scu13, scu16a,scu16b}, Totaro \cite{to08} and Schr\"oer \cite{schr10}, among others\,\e\footnote{\e See \url{https://math.stackexchange.com/questions/4681114}.}\,.

\end{enumerate}
\end{remarks}

\smallskip

Another problem that is discussed in the arithmetic theory of (not necessarily linear) algebraic groups over a number field $k$ is the question of {\it strong approximation}. In \cite{bd}, Borovoi and Demarche studied strong approximation for homogeneous spaces of the form $X=G/\lbe H$, where $G$ is an arbitrary connected algebraic $k$-group and $H$ is a connected $k$-subgroup of $G$. Under certain hypotheses, they essentially showed that the closure of $X(k)$ in the set of adelic points of $X$ is the orthogonal complement, under the Brauer-Manin pairing, of the following subquotient of $\br X$:
\begin{equation}\label{brxg}
\br_{\be a}(\lbe X,G\le)=\krn\be[\e\br X\to \br\e G\to \br\e G_{\ksep}]/\e\img\be[\e\br k\to \br X\le].
\end{equation}
See \cite[Theorem 1.4]{bd}. By work of Demarche \cite[Theorem 4.14]{d11}, if both $H$ and $L=\krn\alb_{\le G}$ \eqref{diag} are linear and {\it reductive}, then \eqref{brxg} is naturally isomorphic to $H^{\le 2}(k, (\be C_{\be X}^{\e\prime}\be)^{\lbe D}\lle)$, where $(\be C_{\be X}^{\e\prime}\be)^{\lbe D}$ is a certain $3\e$-term complex of \'etale sheaves on $k$ which fits into a distinguished triangle
\[
C_{\be X}^{\lle D}\to (\be C_{\be X}^{\e\prime}\be)^{\lbe D}\to{\rm NS}_{A/k}[-1]\to C_{\be X}^{\lle D}[1],
\]
where $C_{\be X}^{\lle D}$ is closely related to a $3\e$-term complex $(\be C_{\be X}^{\e\prime\prime}\lle)^{\lbe D}$ \cite[p.~123, line 1]{d11} of {\it pure Deligne $k$-$1$-motives}. 

\smallskip

What about the function field version of the above question? In contrast with the number field case \cite[\S7.1]{bd}, we can no longer reduce to the case where $L=\krn \alb_{\le G}$ is reductive. Consequently, the $k$-group $Z_{G}$ in \eqref{diag2} will {\it not}, in general, be a semiabelian $k$-variety, and the analog of the complex $(\be C_{\be X}^{\e\prime\prime}\lle)^{\lbe D}$ mentioned above will {\it not} consist of pure Deligne $k$-$1$-motives. This fact was an additional motivation for constructing the category of $k$-$1$-motives described in \cite{ga23}.

\begin{remarks}\label{rmn}\indent
\begin{enumerate}
\item[(a)] Here I stress the fact that the anti\e-\le Chevalley and the Albanese decompositions of a given smooth and connected algebraic $k$-group $G$ are {\it equally} important for investigating the arithmetic of $G$. See, for example, the discussion after diagram \eqref{diag} above. Now, although the category $\s{M}_{k,\e 1}$ of $k$-$1$-motives constructed in \cite{ga23} is built upon Chevalley\e-\le type decompositions of commutative algebraic $k$-groups, anti\e-\le Chevalley decompositions are likely to play a role in establishing duality theorems for the objects of $\s{M}_{k,\e 1}$ in the function field case. I now try to explain why this is so.

If $k$ is any global field, the Kummer theory of Deligne $k$-$1$-motives plays a central role in the proofs of the global duality theorems for such objects given in \cite{hsz1} and \cite{ga09}. Now, if $M=[T^{\lle D}\!\to G\,]\in \s{M}_{k,\e 1}^{\le\rm Del}$ is a Deligne $k$-$1$-motive, then the Kummer theory of $M$ is relatively simple for the following reason: for any positive integer $n$, the multiplication\e-\le by\e-\le $n$ map $[\lle n\lle]$ on both $T^{D}$ and the semiabelian $k$-variety $G$ is an isogeny (i.e., a $k$-morphism with finite kernel and cokernel). It follows that the cohomology of the mapping cone of $[\lle n\lle]\colon M\to M$ is nontrivial only in degree $-1\e$\,\footnote{\e The corresponding fppf sheaf is a finite $k$-group scheme $T_{\Z/\lbe n\le\Z}(M\le)=H^{-1}(M\otimes^{\mathbb L}\Z/\lbe n\le\Z)$ which is an extension of $\,T^{D}\!/n$ by $G_{n}$. See \cite[bottom of p.~208]{ga09}.}\,. Now, the latter is certainly {\it not} the case if $M\in \s{M}_{k,\e 1}$ is an arbitrary (generalized) $k$-$1$-motive and $n$ is divisible by $p$ in the function field case. Thus, it may be quite challenging to establish global duality theorems for an arbitrary $M\in \s{M}_{k,\e 1}$ (in the function field case) via a direct analog of the proof for Deligne $k$-$1$-motives alluded to above. A more sensible approach might be the following: first use the anti\e-\le Chevalley decomposition of the degree $0$ component of $M$, as in the particular case discussed in \cite[Remark 3.7]{ga23}, to obtain an appropriate d\'evissage (in the relevant derived category) of $M$ in terms of Deligne $k$-$1$-motives and  commutative affine algebraic groups (or duals of such), and then appeal to the known global duality theorems for Deligne $k$-$1$-motives \cite{ga09} and commutative affine algebraic groups \cite{ros18} to deduce continuity and perfectness of the appropriate pairings for $M$. Note that, in this alternative approach, the considerations pointed out in Remark \ref{dev} will be of central importance.

\smallskip

\item[(b)] For brevity, a smooth, connected, commutative, unipotent and unirational algebraic $k\e$-group will be called a {\it group of Achet type}. A basic example of such a group (already encountered in Example \ref{ex2} above) is the following: if $\kp\!/k$ is a nontrivial finite and purely inseparable field extension and $T$ is a non\le-\le zero $k$-torus, then the $k$-wound group (apparently, first studied by Oesterl\'e \cite[VI, \S5, pp.~70-73]{oes}):
\begin{equation}\label{ukkt}
O(\kp,T\e)\defeq\rkk(\le T_{\kp}\be)/\le T
\end{equation}
is a group of Achet type \cite[VI, \S5.1, Lemma, p.~70]{oes}. Now, if $U$ is any group of Achet type over a global function field $k$, then the pairing \eqref{sp} exists and is isomorphic to a pairing
\begin{equation}\label{acp}
\Sha^{\le 1}\lbe(\g,U(\ksep)\le)\times \Sha^{\e 1}\lbe(\g,\pic\e U_{\lbe\lle\ksep}\be)\to\Q/\Z
\end{equation}
(cf. Example \ref{ex2}). As noted in Example \ref{ex2}\,, the above pairing is likely to be perfect if $U$ is the $k$-group $O(\kp,T\e)$ \eqref{ukkt}. Now, a natural question is whether or not the pairing \eqref{acp} is perfect for {\it every} group of Achet type $U$. Since $\pic\e U_{\lbe\lle\ksep}$ agrees with ${\rm Ext}^{1}(U_{\lbe\lle\ksep},\G_{m,\e\ksep}\!\lle)$ by \cite[Theorem 2.8(i)]{ach19b}, the answer to the above question may well be affirmative. At any rate, the pairing \eqref{acp} should be investigated for groups of Achet type other than those of the form \eqref{ukkt}, e.g., quotients $O(k^{\lle\prime\prime}\be/\kp,T\le)\defeq O(k^{\lle\prime\prime},T\e)/\le O(k^{\lle\prime},T\e)$, where $k^{\lle\prime\prime}\be/\kp\be/k$ is a tower of nontrivial finite and purely inseparable field extensions.
\end{enumerate}
\end{remarks}
\subsection{Quasi\le-\le reductive groups}\label{qrg} If $k$ is any field, every smooth, connected and affine algebraic $k$-group $G$ is an extension $1\to\s R_{us,\le k}(G\le)\to G\to G^{\lle\rm q\lle r}\to 1$, where $\s R_{us,\le k}(G\le)$ is the $k$-split unipotent radical of $G$ and $G^{\lle\rm qr}\defeq G/\s R_{us,\le k}(G\le)$ is the maximal quasi\le-\le reductive quotient of $G$ \cite[Definition C.2.11, p.~599]{cgp}. It follows that, if $k$ is a global field, then many questions about the arithmetic of $G$ can be reduced to similar questions about the arithmetic of $G^{\lle\rm qr}$\,\e\footnote{\e This is a familiar reduction step in the number field case, where $G^{\le\rm q\lle r}=G^{\le\rm red}$ is, in fact, reductive. See, e.g., \cite[\S7.1]{bd}.}\,.

Now assume that $G$ is quasi\le-\le reductive, i.e., its $k$-unipotent radical $\s R_{u,\le k}(G\le)$ is $k$-wound \cite[Corollary B.3.5, p.~574]{cgp}. Let $G^{\le\rm psr}\defeq G/\s R_{u,\le k}(G\le)$ be the maximal pseudo\le-\le reductive quotient of $G$ and set $G_{1}\defeq G/\s D(\s R_{u,\le k}(G\le))$. Then there exists a canonical exact and commutative diagram of smooth, connected and affine algebraic $k$-groups:
\[
\hskip .75cm\xymatrix{&	\s D(\s R_{u,\le k}(G\le))\ar@{^{(}->}[d]\ar@{=}[r]&\s D(\s R_{u,\le k}(G\le))\ar@{^{(}->}[d]&&\\
1\ar[r]&\s R_{u,\le k}(G\le)\ar@{->>}[d]\ar[r]&G\ar@{->>}[d]\ar[r]&G^{\le\rm psr}\ar@{=}[d]\ar[r]&1\\
1\ar[r]&\s R_{u,\le k}(G\le)/\s D(\s R_{u,\le k}(G\le))\ar[r]&G_{1}\ar[r]&G^{\le\rm psr}\ar[r]&1,
}
\]
where $\s R_{u,\le k}(G\le)/\s D(\s R_{u,\le k}(G\le))=\s R_{u,\le k}(G_{1}\le)$ is commutative and unipotent (but not necessarily $k$-wound, so that $G_{1}$ may not be quasi\le-\le reductive). The middle column of the above diagram expresses $G$ as an extension of $G_{1}$ by the $k$-wound group $\s D(\s R_{u,\le k}(G\le))$, which is either trivial or has a derived length strictly less than that of $\s R_{u,\le k}(G\le)$. On the other hand, the maximal quasi\le-\le reductive quotient $G_{1}^{\le\rm q\lle r}$ of $G_{1}$ is an extension $1\to \s R_{u,\le k}(G_{1}\le)/\s R_{us,\le k}(G_{1}\le)\to G_{1}^{\lle\rm q\lle r}\to G^{\le\rm psr}\to 1$, where the left-hand group is commutative and $k$-wound.

The above discussion motivates the study of the arithmetic of the quasi\le-\le reductive algebraic $k$-groups $G$ which are extensions
\begin{equation}\label{ugp}
1\to U\to G\to P\to 1,
\end{equation}
where $U$ is {\it commutative and $k$-wound} and $P$ is pseudo\le-\le reductive.

An interesting problem is that of computing the Tamagawa number quotient $\tau(G\le)/\tau(U\le)\tau(P\le)$ associated with \eqref{ugp}. If the extension \eqref{ugp} is {\it central}, then this problem was solved in \cite[III, Theorem 5.3, p.~41]{oes} (see also \cite[Remark 1.3.7]{cf}).

Another interesting question is whether or not Rosengarten's formula \cite[Theorem 1.1]{ros20} 
\begin{equation}\label{rf}
\tau(G\le)=\#\le{\rm Ext}^{1}(G,\G_{m,\e k})/\#\Sha^{1}\be(G\le),
\end{equation}
which is valid if $G$ is either commutative or pseudo\le-\le reductive, holds as well for the extension $G$ \eqref{ugp}, at least if this extension is central.

The above questions, which may be difficult to answer in the stated generality, should certainly be investigated for groups of {\it Lourdeaux\le-\le Achet type}. By definition, these are the extensions $G$ \eqref{ugp} of a group of Lourdeaux type $P$ (see subsection \ref{glt}) by a group of Achet type $U$ (see Remark \ref{rmn}(b)). Since, as noted previously, $P$ and $U$ satisfy Sansuc's additivity lemma \cite[Lemma 6.6]{san} (respectively, are unirational), a natural question is whether or not the groups of Lourdeaux\le-\le Achet type have the same properties, at least when \eqref{ugp} is central.

I now describe an interesting class of groups of Lourdeaux\le-\le Achet type.

Let $H$ be an arbitrary reductive algebraic $k$-group and let $Z$ be a (possibly non\le-\le smooth) central $k$-subgroup of $H$. By \cite[Proposition 2.2]{bga14}, the reductive $k$-group $H/Z$ admits a $t$-resolution, i.e., there exists a central extension $1\to T\to E\to H/Z\to 1$, where $T$ is a $k$-torus and $E$ is a reductive $k$-group such that $\s D(E\le)$ is simply connected. Now let $\kp\!/k$ be a nontrivial finite and purely inseparable field extension of degree $p^{\le n}$. Then the given $t$-resolution induces a central \cite[Proposition A.5.15(1), p.~520]{cgp} extension of smooth, connected and affine algebraic $k$-groups
\[
1\to O(\kp,T\le)\to \rkk(\lbe E_{\kp}\be)/\le T\to\rkk(\lbe H_{\kp}\be/Z_{\kp}\be)\to 1,
\]
where $O(\kp,T\le)$ is Oesterl\'e's group \eqref{ukkt}. Let $G$ be the pullback of the above extension along the closed immersion $\rkk(H_{\kp}\be)/\rkk(Z_{\kp}\be)\into \rkk(\lbe H_{\kp}\be/Z_{\kp}\be)$. Then $G$ is a group of Lourdeaux\le-\le Achet type with associated (central) extension
\begin{equation}\label{ogp}
1\to O(\kp,T\le)\to G\to\rkk(H_{\kp}\be)/\rkk(Z_{\kp}\be)\to 1.
\end{equation}
Regarding the above group $G$, the following question seems reasonable\e: if the group $\rkk(H_{\kp}\be)/\rkk(Z_{\kp}\be)$ in \eqref{ogp} satisfies
the function field analog of Theorem \ref{sat}\,, is the same true for the (non\le-\le reductive) quasi\le-\le reductive $k$-group $G$ \eqref{ogp}? 

The simplest version of the above question occurs when $T=\G_{m,\le k}$, on account of the fact that $H^{\lle 1}\lbe(k,O(\kp,\G_{m,\le k}\lle))\simeq\br\lbe(\lbe k\lle)_{p^{\le n}}$ satisfies the Hasse principle.

\begin{remark} In \cite{ros21}, Rosengarten showed that the formula \eqref{rf} fails to hold for non\e-\le commutative unipotent $k$-groups (more precisely, he constructed specific $2$\le-\le dimensional $k$\le-\le wound counterexamples). Further, such groups can exhibit other pathologies \cite[Theorems 1.8 and 1.12]{ros21}. Presumably, the above led him to state in \cite{ros20} that ``commutative and
pseudo\le-\le reductive groups behave nicely, and everything else is pathological". Regarding the latter statement, it would be interesting to determine whether or not the groups of Lourdeaux\e-\le Achet type defined above, which in general are neither commutative nor pseudo\le-\le reductive, can really exhibit pathologies of the type discussed in \cite{ros21}.
\end{remark}

\section{One\le-\le motives and the BSD conjecture}\label{bsd1}

In \cite{blo80}, Bloch gave a volume\e-\le theoretic formulation of the Birch and Swinnerton\e-\le Dyer (BSD) conjecture for an abelian variety $A$ over a number field $k$. More precisely, he showed

\begin{theorem}\label{bt} \cite[Theorem 1.17\e]{blo80} Let $A$ be an abelian variety over a number field $k$ such that
\begin{enumerate}
\item[(i)] The $L$-function $L(A,s)$ of $A$ has a zero at $s=1$ of order $r={\rm rank}_{\le\Z}A(k)$, and
\item[(ii)] $\Sha^{1}\lbe(A\le)$ is finite.
\end{enumerate}
Further, let $X$ be the semiabelian {\rm variation} of $A$ defined in \cite[(0.3)]{blo80}. Then the full {\rm BSD} conjecture for $A$ holds if, and only if, the Tamagawa number formula
\begin{equation}\label{tf}
\tau(X\le)=\frac{\#\le\pic(X\le)_{\rm tors}}{\#\Sha^{1}\lbe(X\le)}
\end{equation}
holds.
\end{theorem}

We begin this section by recasting (a part of\e) Bloch's proof in terms of Ferrand's pinching construction \cite{fer03}. More precisely, we will show that the semiabelian $k$-variety $X$ in the statement of the theorem can be obtained by pinching the dual abelian variety $A^{\lle t}$ at appropriate (closed) points.

\smallskip

First, we briefly describe Ferrand's general pinching construction in the particular case that is relevant here.

Let $A$ be an abelian variety over a global field $k$. Then $A\simeq{\rm Pic}^{\le 0}_{\lbe A^{\lle t}\be/k}$, where $A^{\lle t}\simeq\extk^{1}\be(A,\G_{m,\le k})$ is the dual abelian variety. Let $Z^{\le\prime}$ be a finite closed $k$-subscheme of $A^{\lle t}$, let $Z$ be a finite $k$-scheme and let $\psi\colon Z^{\le\prime}\to Z$ be a schematically dominant  $k$-morphism, i.e., the induced map $\s O_{\lbe Z}\to \psi_{*}\s O_{\lbe Z^{\prime}}$ is injective. By work of Ferrand \cite{fer03}, the data $\Delta=(Z^{\le\prime},Z, \psi)$ determines a pushout diagram in the category of $k$-schemes 
\begin{equation}\label{pin}
\xymatrix{Z^{\le\prime}\,\ar@{^{(}->}[r]\ar[d]^{\psi}&A^{\lle t}\ar[d]^{\varphi}\\
Z\,\ar@{^{(}->}[r]&Y\be(\be\Delta\be),
}
\end{equation}
where $Y\be(\be\Delta\be)$ is proper and geometrically integral, the horizontal maps are closed immersions, $\varphi$ is finite and surjective and the morphism $A^{\lle t}\lbe\setminus Z^{\le\prime} \to Y\be(\be\Delta\be) \setminus Z$ induced by $\varphi$ is an isomorphism. We say that 
\emph{$Y\be(\be\Delta\be)$ is obtained by pinching $A^{\lle t}$ along $Z^{\le\prime}$ via $\psi$}. Now consider the artinian $k$-algebras $B^{\le\prime}=\s O_{\lbe Z^{\prime}}(Z^{\le\prime}\le)$ and $B=\s O_{\lbe Z}(Z\le)$ and set
\begin{equation}\label{lgp}
L_{\lbe\Delta}=R_{B^{\le\prime}\!/k}(\G_{\le m,\le B^{\prime}}\be)/R_{B/k}(\G_{\le m,\le B}).
\end{equation}
Then $L_{\lbe\Delta}$ is a smooth, connected and affine algebraic $k$-group. By work of Brion \cite{bri15}, the pinching diagram \eqref{pin} induces an exact sequence of commutative algebraic $k$-groups
\[
0\to L_{\lbe\Delta}\to{\rm Pic}_{\e Y\be(\lbe\Delta\lbe)\lbe/k}\to {\rm Pic}_{ A^{\lbe t}\be/k}\to  0.
\]
Now, via the isomorphism $A\isoto{\rm Pic}^{\le 0}_{\lbe A^{\lbe t}\be/k}$, the above sequence induces an exact sequence of smooth, connected and commutative algebraic $k$-groups
\begin{equation}\label{deg}
\mathcal E_{\lbe\Delta}\colon 0\to L_{\lbe\Delta}\to G_{\be\Delta}\to A\to 0,
\end{equation}
where $G_{\be\Delta}={\rm Pic}^{\le 0}_{\e Y\be(\lbe\Delta\lbe)\lbe/k}$.

\smallskip

We view the $k$-scheme $Y\be(\be\Delta\be)$ in \eqref{pin} as a {\it singularization}\,\,\footnote{\e Other possible names are {\it singular deformation} and {\it singular modification}, but we disregarded these terms to avoid confusion with other meanings given to the words {\it deformation} and {\it modification} in the literature. For similar reasons, we called $\mathcal E_{\lbe\Delta}$ a ``variation" of $A$ rather than a degeneration of $A$. Incidentally, note that a singularization (as defined above) is both an alteration and a modification in the sense of \cite[2.17 and 2.20]{dj96}. However, in constrast to \cite{dj96}, our aim here is not to resolve singularities, but rather to create some.} of $A^{\lle t}\e$ determined by the data $\Delta=(Z^{\le\prime},Z, \psi)$\,\footnote{\e This singularization is rather ``mild" since the $k$-scheme $Z^{\le\prime}$ in \eqref{pin} is {\it finite} over $k$ and not, e.g., an abelian subvariety of $A^{\lle t}$. For other (more complicated) instances of pinching, see, e.g., \cite{schr07} and \cite{fs20}.} and the extension $\mathcal E_{\lbe\Delta}$ \eqref{deg} as the associated {\it variation} of its connected Picard scheme ${\rm Pic}^{\le 0}_{\lbe A^{\lle t}\be/k}\simeq A$.

\medskip

We now return to the setting of Theorem \ref{bt}\,, i.e., $k$ is a number field and conditions (i) and (ii) of that theorem hold.

\smallskip

If $r={\rm rank}_{\le\Z}\le A(k)=0$, i.e., $A(k)$ is finite, then $A(k)$ is discrete and cocompact in $A(\A_{k})$ (in fact, $A(\lbe\A_{k}\lbe)$ is compact) and the Tamagawa number $\tau(\lbe A)\defeq{\rm vol}(A(\A_{k})/A(k))$ is defined. In this case (i.e., $r=0$), the semiabelian $k$-variety $X$ constructed in \cite[(0.3)]{blo80} equals $A$ and the formula \eqref{tf} with $X=A$ is, indeed, equivalent to the BSD conjecture (cf. \cite[comments after Theorem 1.1]{ros20}). In the case $r>0$, a close examination of the proof of \cite[Theorem 1.17\e]{blo80} reveals that the author has effectively identified a set of free generators of $A(k)$ onto a single point, so as to be in a situation similar to the case $r=0$. Indeed, the resulting semiabelian variation $X\le$ of $A$ has the property that $X(k)$ is discrete and cocompact in $X(\A_{k})$ \cite[Theorem 1.10]{blo80}. The identification alluded to above can be expressed in terms of Ferrand's pinching construction as follows:

\smallskip

Let $k$ again be any global field and let $\mathcal B=\{P_{1},\dots,P_{\lbe r}\}$ be a subset of $A^{\lle t}\lbe(k)$ whose image in $A^{\lle t}(k)/({\rm tors})$ is a $\Z\e$-basis of this quotient. Further, let $O\in A^{\lle t}(k)$ be the neutral element for the group law on $A^{\lle t}$. Then the data
\begin{equation}\label{db}
\Delta(\mathcal B\le)\defeq(\e\textstyle{\coprod_{\e P\e\in\e \{O\}\le\cup\le \mathcal B}}\,\spec\kappa(P),\e\spec k\le, \textstyle{\prod_{\e P\e\in\e \{O\}\le\cup\le \mathcal B}}\,{\rm Id}\le),
\end{equation}
where ${\rm Id}\colon \spec\kappa(P)\to \spec k$ is the identity morphism for each $P\e\in\e \{O\}\le\cup\le \mathcal B$, determines a pinching diagram
\begin{equation}\label{pin2}
\xymatrix{\displaystyle{\coprod_{\le P\e\in\e \mathcal M^{\lle t}}}\spec\kappa(P)\,\ar@{^{(}->}[r]\ar[d]&A^{\lle t}\ar[d]\\
	\spec k\,\ar@{^{(}->}[r]&Y\lbe(\lbe\mathcal B\le),
}
\end{equation}
where $Y\lbe(\lbe\mathcal B\le)\defeq Y\lbe(\lbe\Delta(\lbe\mathcal B\le)\be)$ and the top horizontal map is the canonical closed immersion. The above is one of the simplest types of pinching constructions. In the case ${\rm dim}\e A^{\lle t}=1$ and $r>0$, $Y\lbe(\lbe\mathcal B\le)$ is a (singular) projective curve with $r$ loops (or ``petals") attached to a single point. Now, the $k$-group $L_{\lbe\Delta(\mathcal B\le)}$ \eqref{lgp} determined by the data $\Delta(\mathcal B\le)$ \eqref{db} can be identified with the split $k$-torus $T=\G_{m,\le k}^{\le r}$ and the sequence $\mathcal E_{\lbe\Delta(\mathcal B\le)}$ \eqref{deg} can be identified (in the number field case) with Bloch's sequence $0\to T\to X\to A\to 0$ \cite[(0.3)]{blo80}, i.e., the $k$-group $X$ in \eqref{tf} is isomorphic to $G_{\be\Delta(\mathcal B\le)}={\rm Pic}^{\le 0}_{\e Y\lbe(\lbe\mathcal B\le)\be/k}$.

\begin{remark} It is stated in \cite{fer03} that the ideas contained therein were developed in the early 1970's but were only published in [loc.cit.], i.e., about three decades later (!). Thus, Bloch may well have been unaware of Ferrand's pinching construction when he wrote \cite{blo80}.
\end{remark}

I will not discuss here other parts of the proof of \cite[Theorem 1.17\e]{blo80} (apart from the brief comments contained in the last section), since those parts are not directly relevant to this letter.

I will now show how Ferrand's pinching construction and $k$-$1$-motives (more precisely, the pure $k$-$1$-motives that arise from \eqref{deg}) can be combined to produce a new method for studying the structure of $\Sha^{\le 1}\be(A\le)$ that differs significantly from the classical descent method.

\smallskip

Let $\Delta=(Z^{\le\prime},Z, \psi)$, $G_{\be\Delta}$ and $\mathcal E_{\lbe\Delta}$ be as above and consider the (smooth) $k$-$1$-motive $M_{\Delta}$ which is dual to $(G_{\be\Delta},\mathcal E_{\Delta})\e$, i.e., $M_{\Delta}$ is the $k$-$1$-motive
\begin{equation}\label{md}
M_{\Delta}=(G_{\be\Delta},\mathcal E_{\Delta})^{\vee}=\left[L_{\be \Delta}^{\! D}\overset{\!v_{\lbe\Delta}}{\lra}A^{\lle t}\le\right].
\end{equation}
The map $v_{\Delta}$ above induces a map $v_{\Delta}^{(1)}\colon H^{1}(k,L_{\be \Delta}^{\! D})\to H^{1}(k,A^{\lle t})$. Since, by \cite[Theorem 3.11(2)]{how}, the preceding construction behaves well under arbitrary extensions of $k$ (specifically, under the extensions $k_{v}/k$ for every place $v$ of $k$), we conclude that $v_{\Delta}^{(1)}$ induces a homomorphism of abelian groups 
\begin{equation}\label{sha}
\Sha(v_{\Delta}^{(1)})\colon \Sha^{\le 1}\be(L_{\be \Delta}^{\! D})\to \Sha^{\le 1}\lbe(\lbe A^{\lle t}\le).
\end{equation}
The above map, whose source is finite by \cite[Proposition 5.16.1]{ros18}, can be described explicitly (in terms of torsors, for example) by unwinding the various definitions involved in its construction. 

\smallskip

We conclude that:

\smallskip

{\it Every singularization of $A^{\lle t}$ of the form $Y(\be\Delta\lbe)$ \eqref{pin} defines a finite subgroup of $\Sha^{\le 1}\lbe(\lbe A^{\lle t}\le)$, namely $\img\le \Sha(v_{\Delta}^{(1)})$, where $\Sha(v_{\Delta}^{(1)})$ is the map \eqref{sha}.}

\smallskip

\begin{remark} Note that the map \eqref{sha} induced by the singularization $Y\lbe(\lbe\mathcal B\le)$ \eqref{pin2} which produces Bloch's semiabelian $k$-variety $X$ is the {\it zero map}, since $\Sha^{\le 1}\be(L_{\be \Delta}^{\! D})\simeq \Sha^{\le 2}\lbe(\G_{m,\le k}^{\le r})^{*}=0$ by the Brauer-Hasse-Noether theorem. 
\end{remark}

\begin{example} \label{rr} Let $n\geq 1$ be an integer and set $\Delta=\Delta(n)=(A^{\lle t}_{\lle\lle n},\spec k, \psi)$, where $A^{\lle t}_{\lle\lle n}$ is the $n$-torsion subgroup of $A^{\lle t}$ and $\psi\colon A^{\lle t}_{\lle\lle n}\to\spec k$ is its structural morphism. Set $L_{\le(n)}=L_{\Delta(n)}= R_{A^{\lle t}_{\lle n}/k}(\G_{m,\e A^{\lbe t}_{\lle n}}\be)/\G_{m,\le k}$. Then the B\'egueri embedding $A_{\le n}\into R_{A^{\lle t}_{\lle n}/k}(\G_{m,\e A^{\lbe t}_{\lle n}}\be)$ \cite{ha20} induces a $k$-morphism $A_{\le n}\to L_{(n)}$ which (at least if $n$ is prime to $\car k$) induces a map $L_{(n)}^{\be D}\to A^{\lle t}_{\lle\lle n}$. The latter map induces, in turn, a map $\theta_{\lle n}\colon \Sha^{1}(L_{\lbe(n)}^{\! D})\to 
\Sha^{1}\lbe(A^{\lle t}_{\lle\lle n})$. It can be checked that $\Sha(v_{\!\Delta(n)}^{(1)})$ \eqref{sha} is the composition
\begin{equation}\label{comp}
\hskip .5cm \Sha^{1}(L_{\lbe(n)}^{\! D})\overset{\!\!\theta_{\lbe n}}{\lra} \Sha^{1}\lbe(A^{\lle t}_{\lle\lle n})\into {\rm Sel}^{(n)}(A^{\lle t})\onto \Sha^{1}\lbe(A^{\lle t}\le)_{n}\into \Sha^{1}\lbe(A^{\lle t}\le),
\end{equation}
where
\begin{equation}\label{sel}
{\rm Sel}^{(n)}(\lbe A^{\lle t}\lle)=\krn\!\!\be\left[\lle H^{1}(k,A^{\lle t}_{\lle\lle n})\to\prod_{{\rm all}\,v}H^{1}(k_{\le v},A^{\lle t})\right]\!\be
\end{equation}
is the $n$-Selmer group of $A^{\lle t}$. In general, the Selmer group \eqref{sel} is difficult to compute (see, e.g., \cite{dss00}).
By constrast, the group $\Sha^{1}(L_{\lbe(n)}^{\! D})$ can often be computed effectively. For example, by recent work of Borovoi and Kaletha \cite{bk}, $\Sha^{1}(L_{\lbe(n)}^{\! D})$ is effectively computable in many cases if $L_{\lbe(n)}=T$ is a $k$-torus (e.g., $A^{\lle t}_{\lle\lle n}$ is \'etale over $k$ \cite[Proposition 4.10]{bri15}). In this case
$\Sha^{1}(L_{\lbe(n)}^{\! D})=\Sha^{1}(T^{\lle D}\lle)\simeq\Sha^{\le 2}\lbe(T\le)^{*}$ and, by \cite[Proposition 7.3.6(2)]{bk}, the following holds: if  $X_{*}(T\le)$ and $F\subset k(A^{\lle t}_{\lle\lle n})$ are, respectively, the cocharacter lattice and minimal Galois splitting field of $T$, then  
\[
\Sha^{\le 2}\lbe(T\le)=\cok\!\!\left[\e\bigoplus_{{\rm all}\,v}X_{*}(T\le)_{\le\Gamma(v),\e{\rm tors}}\overset{\prod\be{\rm Cor}}{\lra} X_{*}(T\le)_{\le\Gamma,\e{\rm tors}}\right]\!\!,
\]
where, for each place $v$ of $k$, $\Gamma(v)\subset \Gamma={\rm Gal}(F/k)$ is the decomposition group of $v$ in $F/k$. Thus, if the map $\theta_{\lbe n}$ in \eqref{comp} is nonzero, then the above argument yields a new method for computing nontrivial elements in $\Sha^{1}\lbe(A^{\lle t}_{\lle\lle n})$ and therefore in ${\rm Sel}^{(n)}(A^{\lle t})$ (and also in $\Sha^{1}(A^{\lle t}\le)_{n}$ if the composite map $\Sha^{1}(L_{\lbe(n)}^{\! D})\overset{\!\!\theta_{\lbe n}}{\lra}  {\rm Sel}^{(n)}(A^{\lle t})\onto \Sha^{1}(A^{\lle t}\le)_{n}$ is nonzero).
\end{example}

\section{Concluding remarks}

I believe to have shown that the new theory of $k$-$1$-motives is {\it necessary} for developing a duality theory for commutative algebraic that works equally well in the number field and function field cases; that it is not ``generalization for its own sake" and that it does really have (multiple) potential applications. Now, regarding the contents of section \ref{bsd1}\e, I would like to add the following comments. Firstly, I have shown that, in certain contexts, {\it adding} (rather than removing) singularities can be fruitful. For another illustration of this fact, see, e.g., \cite{deo18}. Secondly, as I understand it, the classical descent method arises from an attempt to study the arithmetic of an abelian variety $A$ by studying the arithmetic of certain {\it affine} algebraic groups associated to $A$, namely the $n$\le-\le torsion subgroups $A_{\le n}$ of $A$, which are necessarily {\it finite} over $k$. Certainly, an abelian variety cannot contain other types of affine subgroups. However, I have shown (at least, in what concerns the Tate\le-\le Shafarevich group) that, by creating certain types of singularities on $A$ and using the theory of $k$-$1$-motives, we can study the arithmetic of $A$ by studying the arithmetic of certain other (i.e., non\le-\le finite) types of affine algebraic $k$-groups $L_{\lbe\Delta}$ that are naturally associated to $A$. Further, we can do this for many possible choices of pinching locus $\Delta$, not just the choice $\Delta=\Delta(n)$ discussed in Example \ref{rr}\,. Moreover, this pinching process can be iterated, e.g., we can pinch $A$ at $\coprod_{\e\ell\text{\e\le\lle (prime)\e}\leq\e x}\coprod_{N\le\leq\e y}\coprod_{\e m=0}^{\e N}(\Delta(\ell^{\le m+1})\setminus \Delta(\ell^{\le m}))$ for increasing values of $x, y\in\R$, etc. Now, it seems reasonable to expect that new insights into the BSD conjecture will be gained by pinching Bloch's $k$-scheme $X$ further. I should note, however, that the discussion contained in section \ref{bsd1} is unsatisfactory since it concerns only the abelian variety $A$ over the global field $k$. A fuller discussion of the BSD conjecture certainly involves choosing $S$-models of $A$ and of its variations, where $S$ is the relevant Dedekind base. Bloch worked with the N\'eron $S$-model of $A$ and the N\'eron-Raynaud $S$-model of its variation $X$. This seems like a natural choice and sufficed for his purposes. However, I am not convinced that the N\'eron-Raynaud model will always be the correct $S$-model to choose. The reason is that the process that produces the N\'eron $S$-model of $A$ (namely blowups, group smoothing, etc) is concerned with extending to $S$ the smoothness and group structure of $A$ (and certain other properties), but {\it not} the Albanese\le-\le Picard duality satisfied by $A$. Since the latter plays a central role in the theory of $k$-$1$-motives, it is currently unclear (to me) how exactly to develop a ``correct" theory of $S$-models of pure $k$-$1$-motives. Developing such a theory seems necessary, since any additional pinching that is performed on Bloch's $k$-scheme $X$ immediately raises the question of choosing adequate (possibly non\le-\le regular and non\le-\le separated) $S$-models of the resulting variations of $X$, which induce corresponding $S$-models of the associated dual $k$-$1$-motives \eqref{md}. Incidentally, note that any $S$-model of the map $v_{\Delta}$ in \eqref{md} yields a new approach (i.e., one that does not involve the Selmer group or any of its variants) for studying the growth of the Mordell\le-\le Weil rank of $A^{\lle t}$, as such an $S$-morphism induces (under certain conditions) a map from the group of characters of a group of relative units of a finite extension of $k$ to $A^{\lle t}(k)$.

Problems of the above sort, and related ones, will keep me busy for quite a long time. Consequently, I will not be able to work on any of the problems pointed out in section \ref{aggf}\,. Since I believe that these problems are interesting and important and should be investigated, I hereby offer my help (as a remote co\le-\le advisor, for example) to any doctoral student who may be interested in solving some of these problems.
\end{spacing}


\begin{thebibliography}{[19]}
	
	
\bibitem[Achet19a]{ach19a} Achet, R.:\emph{ Picard group of unipotent groups, restricted Picard functor}. \url{https://hal.science/hal-02063586/document}

\bibitem[Achet19b]{ach19b} Achet, R.:\emph{ Unirational algebraic groups}. \url{https://hal.science/hal-02358528/document}

\bibitem[Blo80]{blo80} Bloch, S.:\emph{ A Note on Height Pairings, Tamagawa Numbers, and the Birch and Swinnerton-Dyer Conjecture}. Invent. Math. {\textbf{58}} (1980), 65--76.

\bibitem[Bo]{bo} Borel, A.:\emph{ Linear algebraic groups}, Second enlarged edition, Springer-Verlag, New York, 1991.


\bibitem[Bor\e 98]{bor} Borovoi, M.:\emph{ Abelian Galois cohomology of
reductive groups.} Mem. Amer. Math. Soc. {\textbf{132}} (1998), no. 626.



\bibitem[BD]{bd} Borovoi, M. and  Demarche, C.:\emph{ Manin obstruction to strong approximation for homogeneous spaces.} Comment. Math. Helv. {\textbf{88}} (2013), 1--54.


\bibitem[BCTS]{bcts} Borovoi, M., Colliot-Th\'el\`ene, J.-L. and Skorobogatov, A.:
\emph{The elementary obstruction and homogeneous spaces}. Duke Math. J. \textbf{141}(2): 321--364.

\bibitem[BGA14]{bga14} Borovoi, M. and Gonz\'alez-Avil\'es, C.D.:
\emph{The algebraic fundamental group of a reductive group scheme over an arbitrary base scheme}. Cent. Eur. J. Math. \textbf{12}, no.4 (2014), 545--558.


\bibitem[BvH]{bvh} Borovoi, M. and van Hamel, J.: Extended
Picard complexes and linear algebraic groups.  J. reine angew.
Math. {\bf{627}} (2009), 53--82.

\bibitem[BK23]{bk} Borovoi, M. and Kaletha, T.: \emph{ Galois cohomology of reductive groups over global fields.} \url{https://arxiv.org/abs/2303.04120}





\bibitem[Bri15]{bri15} Brion, M.:\emph{ Which algebraic groups are Picard varieties?} Sci. China Math. \textbf{58} (2015), no. 3, 461--478.

\bibitem[Bri17a]{bri17a} Brion, M.: \emph{Some structure theorems for algebraic groups}. Algebraic groups: structure and actions, 53--126, Proc. Sympos. Pure Math., \textbf{94}, Amer. Math. Soc., Providence, RI, 2017.




\bibitem[Che]{che} Chernousov, V.:\emph{ The Hasse principle for groups of type $E_{8}$.} Soviet Math. Dokl. \textbf{39} (1989), 592--596.


 
 
\bibitem[CP]{cp} Conrad, B. and Prasad, G.: Classification of pseudo-reductive groups. Annals of Mathematics Studies \textbf{191}, Princeton Univ. Press, 2015.


\bibitem[CGP]{cgp} Conrad, B., Gabber, O. and Prasad, G.: Pseudo-reductive groups. Second Ed. New Math. Monograps \textbf{26}, Cambridge U. Press, 2015.

 
\bibitem[CF]{cf} Conrad, B., \emph{Finiteness theorems for algebraic groups over function fields}. Compos. Math. \textbf{148}, no. 2 (2012), 555--639.
 
 
\bibitem[dJ96]{dj96} de Jong, J.A.:\emph{ Smoothness, semi\le-\le stability and alterations.} Publ. Math. IHES \textbf{83} (1996), 51--93.

\bibitem[D11]{d11} Demarche, C.:\emph{ Une formule pour le groupe de Brauer alg\'ebrique d'un torseur.} J. Algebra \textbf{347} (2011), 96--132.

\bibitem[Dem11]{dem11} Demarche, C.:\emph{ Suites de Poitou-Tate pour les complexes de tore \`a deux termes”}. Internat. Math. Res. Notices (2011), no. 1, p. 135--174 (abridged  version of ``Th\'eor\`emes de dualit\'e pour les complexes de tores”, available at \url{https://webusers.imj-prg.fr/~cyril.demarche/}.



\bibitem[DH20]{dh20} Demarche, C. and Harari, D.:\emph{ Duality for complexes of tori over a global field of positive characteristic}. J. de l'\'Ecole Polytechnique (Math\'ematiques) \textbf{7}, (2020) 831--870.

\bibitem[DH22]{dh22} Demarche, C. and Harari, D.:\emph{
Local-global principles for homogeneous spaces of reductive
groups over global function fields.} Annales H.Lebesgue \textbf{5} (2022), 1111--1149.


\bibitem[DG]{dg}  Demazure, M. and Gabriel, P.:\emph{ Groupes alg\'ebriques}. Tome I: G\'om\'etrie alg\'ebrique, g\'en\'eralit\'es, groupes commutatifs. Avec un appendice Corps de classes local par Michiel Hazewinkel. Masson  \& Cie, \'Editeur, Paris; North-Holland Publishing Co., Amsterdam, 1970.


\bibitem[$\text{SGA3}_{\le\text{new}}$]{sga3}  Demazure, M. and
Grothendieck, A. (Eds.): Sch\'emas en groupes. S\'eminaire de
G\'eom\'etrie Alg\'ebrique du Bois Marie 1962-64 (SGA 3). Augmented and
corrected 2008-2011 re-edition of the original by P.Gille and P.Polo.
Available at \url{http://www.math.jussieu.fr/~polo/SGA3}. Reviewed
at \url{http://www.jmilne.org/math/xnotes/SGA3r.pdf}.

\bibitem[Deo18]{deo18} Deopurkar, A. \emph{ The canonical syzygy conjecture for ribbons.} Math. Z. \textbf{288} (2018), 1157--1164.


\bibitem[DSS00]{dss00} Djabri, Z., Schaefer, E. and Smart, N.:\emph{ Computing the $p$-Selmer group of an elliptic curve.} Trans. Amer. Math. Soc. {\textbf{352}}, no. 12 (2000), 5583--5597.


\bibitem[FS20]{fs20} Fanelly, A. and Schr\"oer, S.:\emph{ Del Pezzo surfaces and Mori fiber spaces in positive characteristic.} Trans. Amer. Math. Soc. {\textbf{373}}, no. 3 (2020), 1775--1843.


\bibitem[Fer03]{fer03}  Ferrand, D.:\emph{ Conducteur, descente et pincement}. Bull. Soc. Math. France \textbf{131} (2003), no. 4, 553--585. 


\bibitem[GA09]{ga09} Gonz\'alez-Avil\'es, C.D.:
\emph{ Arithmetic duality theorems for $1$-motives over function fields.} J. reine angew. Math. \textbf{632} (2009), 203-231.


\bibitem[GA12]{gaq} Gonz\'alez-Avil\'es, C.D.:
\emph{ Quasi-abelian crossed modules and nonabelian cohomology.} Journal of Algebra \textbf{369} (2012), 235--255.


\bibitem[LDAG]{ga23} Gonz\'alez-Avil\'es, C.D.:
\emph{ Local duality theorems for commutative algebraic groups.} Available at
\url{https://arxiv.org/abs/2305.08699}.


\bibitem[Ha20]{ha20} Haine, P. \emph{ The B\'egueri resolution}. \url{https://math.berkeley.edu/~phaine/#research}


\bibitem[HSz05]{hsz1}  Harari, D. and Szamuely, T.:\emph{ Arithmetic
duality theorems for $1$-motives.} J. reine angew. Math. {\bf{578}},
pp. 93-128 (2005), and \emph{Errata}: available from
http://www.renyi.hu/$\sim$szamuely.


\bibitem[HSz08]{hsz2} Harari, D. and Szamuely, T.:\emph{ Local-global principles
for 1-motives}. Duke Math. J. {\bf 143}(3) (2008), 531-557.



\bibitem[How12]{how} Howe, S.: \emph{ Higher genus counterexamples to relative Manin-Mumford.} Master's degree thesis directed by Bas Edixhoven, Universiteit Leiden and Universit\'e Paris-Sud 11, 2012.

\bibitem[Lour]{lord} Lourdeaux, A. \emph{ On geometry of pseudo-reductive groups}. \url{https://arxiv.org/abs/2011.14001}




\bibitem[Miz73]{miz73} Mizutani, H.:\emph{ Hironaka's additive group schemes.}
Nagoya Math. J. \textbf{52} (1973), 85--95.


\bibitem[Oda73]{oda73} Oda, T.:\emph{ Hironaka's additive group scheme.} Number theory, algebraic geometry and commutative algebra in honor of Y.Akizuki, Kinokuniya, Tokyo, 1973, 181--219.



\bibitem[Oes]{oes} Oesterl\'e, J.:\emph{ Nombres de Tamagawa et groupes unipotents en caract\'eristique $p$}. Invent. Math. \textbf{78}  (1984), no. 1, 13--88.


\bibitem[Ros18]{ros18} Rosengarten, Z.:\emph{ Tate duality in positive dimension over function fields}. \url{https://arxiv.org/abs/1805.00522}.

\bibitem[Ros20]{ros20} Rosengarten, Z.:\emph{ Tamagawa numbers and other invariants of pseudo-reductive groups over global function fields}. \url{https://arxiv.org/abs/1806.10723}.

\bibitem[Ros21]{ros21} Rosengarten, Z.:\emph{ Pathological behavior of arithmetic invariant of unipotent groups}. \url{https://arxiv.org/abs/1907.06225}.

 
\bibitem[Ros18]{ros18} Rosengarten, Z.:\emph{ Tate duality in positive dimension over function fields}. \url{https://arxiv.org/abs/1805.00522}.

\bibitem[San]{san} Sansuc, J.-J.:
\emph{ Groupe de Brauer et arithm\'etique des groupes alg\'ebriques
lin\'eaires sur un corps de nombres.}  J. reine angew. Math.
{\textbf{327}} (1981), 12-80.



\bibitem[Schr07]{schr07} Schr\"oer, S.:\emph{ Weak del Pezzo surfaces with irregularity.} Tohoku Math. J. {\textbf{59}} (2007), 293--322.

\bibitem[Schr10]{schr10} Schr\"oer, S.:\emph{ On fibrations whose geometric fibers are nonreduced.} Nagoya Math. J. {\textbf{200}} (2010), 35--57.



\bibitem[Scu13]{scu13} Scully, S.:\emph{ Rational maps between quasilinear hypersurfaces.} Compos. Math. {\textbf{149}}, no. 3, (2013), 333--355.


\bibitem[Scu16a]{scu16a} Scully, S.:\emph{ On the splitting of quasilinear $p\,$-forms.} J. Reine Angew. Math. {\textbf{713}} (2016), 49--83.

\bibitem[Scu16b]{scu16b} Scully, S.:\emph{ Hoffmann's conjecture for
totally singular forms of prime degree.} Algebra and Number Theory \textbf{10}, no. 5 (2016), 1091--1132.


\bibitem[Sko]{sko} Skorobogatov, A.: Torsors and rational points. 2001.

\bibitem[SP]{sp} The Stacks Project. \url{http://stacks.math.columbia.edu}


\bibitem[To08]{to08} Totaro, B.:\emph{ Birational geometry of quadrics in characteristic $2$.} J. Algebraic Geom. \textbf{17}, no. 3, (2008), 577--597.


\bibitem[To13]{to13} Totaro, B.: \emph{Pseudo-abelian varieties}. Ann. Sci. \'Ec. Norm. Sup\'er. (4) \textbf{46} (2013), no. 5, 693--721. 



\bibitem[vH04]{vhd} van Hamel, J. \emph{ Lichtenbaum-Tate duality for varieties over $p\e$-adic fields}. J. Reine Angew. Math. \textbf{575} (2004), 101--134.

\end{thebibliography}
\end{document}